\documentclass[12pt]{amsart}
\usepackage{times}
\usepackage{amscd}
\usepackage[latin1]{inputenc}
\usepackage{latexsym}
\usepackage{amsmath,amssymb,amsthm,amsfonts}
\usepackage[dvips]{graphicx}
\usepackage{epsf}
\newtheorem{teor}{Theorem}[section]
\newtheorem{lemma}[teor]{Lemma}
\newtheorem{propos}[teor]{Proposition}
\newtheorem{corol}[teor]{Corollary}
\theoremstyle{definition}
\newtheorem{defin}[teor]{Definition}
\newtheorem{osserv}[teor]{Remark}

\newcommand{\Q}{\mbox{ $ \mathbb{ Q}$}}

\newcommand{\C}{\mbox{ $ \mathbb{C} $ }}

\newcommand{\F}{\mbox{$\mathcal{F}$}}

\newcommand{\M}{\mbox{ $ M$ }}

\newcommand{\p}{\mbox{$p$}}

\newcommand{\fu}{\mbox{ $ f:M \longrightarrow M $ }}
\newcommand{\id}{\mbox{  \hbox{Ind}}}
\newcommand{\Fl}{\mbox{  $ \mathcal{F} $ }}

\theoremstyle{remark}  %\newtheorem{nota}[teo]{Remark}
\numberwithin{equation}{section}
\begin{document}
\title[Holomorphic dynamics near germs of singular curves]{Holomorphic dynamics
near germs of singular curves}
\author[Francesco Degli Innocenti]{Francesco Degli Innocenti}
\address{Dipartimento di Matematica, Universit\`a di Pisa, via Buonarroti 2, 56127, Pisa.}
\email{degliinno@dm.unipi.it}
\begin{abstract}
Let $M$ be a two dimensional complex manifold, $p \in M $  and \Fl a germ of
holomorphic foliation of \M at $p$. Let $S\subset M$ be a germ of an
irreducible, possibly singular, curve at $p$ in  $M$ which is a separatrix for \Fl.
We prove that if the Camacho-Sad-Suwa index  $\id(\F,S,p)\not \in \Q^+\cup \{0\} $ then there exists
another separatrix for \Fl at $p$. A similar  result is
proved for the existence of parabolic curves for
germs of holomorphic diffeomorphisms near a curve of fixed points.
\end{abstract}

%\title[dynamics near curves of fixed points]{The dynamics of holomorphic
%maps near  curves of fixed points.}
%\author[Filippo Bracci]{Filippo Bracci$^\dag$}
%\thanks{\rm $^\dag$ Partially supported by Progetto MURST di
%Rilevante Interesse Nazionale {\it Propriet\`a geometriche delle
%variet\`a reali e complesse.}} \address{Dipartimento di
%Matematica, Universit\`a di Roma  ``Tor Vergata'', Via della
%Ricerca Scientifica 00133 Roma, Italy.}
%\email{fbracci@mat.uniroma2.it}
 %\subjclass{Primary 32H50, 37F99,
%32S45, 32S65.}
\maketitle
\section{Introduction}
Let $M$ be  a two dimensional complex manifold and \Fl a germ of
holomorphic foliation on $M$ near $p$. In local coordinates the foliation
can be described by the vector field:
$$ A(x,y)\frac{\partial}{\partial x}+B(x,y)\frac{\partial}{\partial y},$$
with $A, B$  suitable holomorphic functions.
A separatrix for \Fl is a non constant holomorphic solution of the system:
$$ \left \{ \begin{aligned} &\dot{x}=A(x,y)\\
                            &\dot{y}=B(x,y)
            \end{aligned}\hspace{2cm} x(0)=y(0)=0.  \right.$$
Obviously the interesting case is when $(0,0)$ is a singularity
for \F.
In the singular case, in the well known paper \cite{CS}, Camacho and Sad proved that there
always exists (at least) one (possibly singular) irreducible separatrix - say $S$ - for \Fl at $(0,0)$.
A natural question is whether the knowledge of this separatrix $S$ allows to infer the existence of
 another separatrix.
There are essentially two types of results, one of local and the other of global flavour.
The first kind of result is essentially a re-formulation of Camacho-Sad theorem
(see the paper by J. Cano \cite{CA}) which says that if $S$ is non singular and
$ \id(\F,S,p)\not \in \Q^{+}\cup \{0\}$ (where $\id(\F,S,p)$ is the index introduced in \cite{CS}) then
there exists another separatrix through $p$. The second type of result requires global
conditions on $S$, like $S$ compact (but possibly singular), globally and locally irreducible and
$S \cdot S < 0$ to provide the existence of another separatrix at some point of $S$ (see the paper by
Sebastiani \cite{SB}).

One aim of this paper is to prove a result of local nature when $S$ is possibly singular, using the
index defined by Suwa \cite{SW}.
We prove:
\begin{teor} \label{teor:miomio}
Let  \M  be a complex two dimensional manifold,  \Fl  a holomorphic  foliation on same open subset of
 $M$, $S \subset M $ a possibly singular curve locally irreducible at a point $p\in M $, such that it is
 a separatrix for   \Fl at $ p $.
 If $ \id(\F,S,p)\not \in \Q^+\cup\{0\} $
then there exists (at least) another separatrix for \Fl at  $p $.
\end{teor}

Abate, Bracci and Tovena \cite{ABT}, \cite{BR}, \cite{BT}   have
recently shown how  to translate results about foliations to holomorphic diffeomorphisms with curves
of fixed points.
The proof of Theorem~\ref{teor:miomio} respects their dictionary and so the results about
the existence of separatrices for foliations can be translated in
results about the existence of parabolic curves for
diffeomorphisms. Using notations of \cite{ABT}, \cite{BR} and
\cite{BT}  we obtain:

\begin{teor}\label{teor:diff}
Let $M$ be a two dimensional complex manifold, \fu a holomorphic
map  such that $\hbox{Fix}(f)=S$ with $S$ a locally
irreducible, possibly singular curve at a point $p\in M$. Assume that $f$ is tangential
on $S$ and  $\id(f,S,p)\not \in \Q^+\cup \{0\} $. Then there exists (at least) a parabolic curve for $f$
at $p$.
\end{teor}
Theorem \ref{teor:diff} has been proved by Abate \cite{AB} in case $S$ is non singular and by  Bracci in \cite{BR} in case $S$ is a
generalized cusp, i.e. of the form $\{x^m=y^n\}$.

I want to sincerely thank professor Filippo Bracci without whose help this work would not have came to be.

\section{Preliminary results}
First of all we have to recall some basic notions about C.S.S. (Camacho-Sad-Suwa) index.
This index was first introduced by Camacho and Sad in \cite{CS} for a complex one codimension singular foliation defined in a neighborhood
of a non singular compact curve  embedded  in a two dimensional complex manifold.
Later Suwa \cite{SW} generalized to a generic possibly singular compact  invariant curve.
The most interesting property of this is the following Index Theorem, that relates the dynamics of \Fl near a curve $S$
to the self intersection number of $S$.
\begin{teor}[Index Theorem]
 Let  $S$ be  a compact curve in a two dimensional complex manifold $M$ invariant by a possibly singular
foliation \Fl,then for every point $p\in S$ there exists a complex  number $\id(\Fl,S,p)\in \C$ depending only on the local behaviour of
\Fl and $S$ near $p$ such that:
\begin{equation*}
\sum_{p\in S}\id(\mathcal{F},S,p)= S \cdot S.
\end{equation*}
\end{teor}
We now recall the behaviour of this index under blow-up.

\begin{propos}\label{propos:blow-up}
Let $M $ be a two dimensional complex manifold, \Fl an holomorphic foliation, $S$ an \Fl-separatrix and $p\in S $ a singularity of $S$.
We indicate by $ \pi: \tilde{M}\longrightarrow M $ the blow-up of $M$ in $p$, by $ \tilde{\mathcal{F}} $ the saturated foliation and
by  $ D :=\pi^{-1}(p) $ and $ \hat{S}:=\overline{\pi^{-1}(S\setminus \{p\})} $ respectively the exceptional divisor and the strict
transform of $S$.
Then $\hat{S} $ is an $\tilde{\mathcal{F}}$ separatrix.
Moreover if $ \{\tilde{p}\}:=D\cap\hat{S} $ then
$$ \hbox{Ind}(\tilde{\mathcal{F}},\hat{S},\tilde{p})=\hbox{Ind}(\mathcal{F},S,p)-m^2 $$
where $ m \geq 1 $ is the multiplicity of $S$ in $p$.
\end{propos}

Cano in \cite{CA} gives an algorithmic proof of Camacho-Sad result introducing a particular class of points that we will often  use.
\begin{defin}
Let $M$ be a two dimensional complex manifold, \Fl an holomorphic foliation and $S$ a local separatrix for \Fl.
\begin{itemize}
\item
We say that a point $p\in S$ is of type $(C_1)$ if $S$ is nonsingular at $p$ and
$$ \id(\Fl,S,p)\not \in \Q^+ \cup \{0\}.$$
\item
We say that a point $p\in S$ is of type $(C_2)$ if $S$ has two nonsingular branches $S_0$, $S_1$ at $p$, intersecting trasversally at
$p$, and there exists a real number $r>0$ such that
$$ \begin{aligned}
&\id(\mathcal{F},S_0,p)\not \in \Q_{\geq -\frac{1}{r}}=\{a\in \Q : a\geq -\frac{1}{r}\}\\
&\id(\mathcal{F},S_1,p)\in \Q_{\leq -r}=\{a\in \Q : a\leq -r\}.
\end{aligned}$$
\end{itemize}
\end{defin}

According to Definition 7.6 of \cite{BR} and \cite{CA} we have:
\begin{defin}
A point $p\in S \subset M$ where  $S$ is an \Fl-invariant curve is said to be an \textbf{appropriate singularity} for \Fl if after
a finite number of blow-ups there exists a $(C_1) $ or $(C_2)$ point on the total transform.
\end{defin}
The importance of this class of points is given by the following result:
\begin{propos}[\cite{BR}, \cite{CA}]
If $p\in S \subset M$ is an appropriate singularity for a foliation \Fl, then at least another separatrix trough $p$ for \Fl exists.
\end{propos}

\section{Proof of the result}
In order to get Theorem \ref{teor:miomio} we will concentrate our attention on the  particular class of points introduced in the previous
section.
The upshot is  to show that under the hypotheses of Theorem \ref{teor:miomio} the point $p$ is an appropriate singularity.

We know that the resolution of curves singularities theorem \cite{LAU} ensures that after a finite
 number of blow-ups we have the geometric structure required for the existence of $(C_1)$ or $(C_2)$ points.
To conclude we have to analize the C.S.S. index under  this process.
The behaviour of the index is strongly related to the  evolution of the geometric structure under blow-up.
We can divide the proof in two steps:
\begin{enumerate}
\item study of the geometric structure under the resolution of singularities,
\item study of the C.S.S. index under this process.
\end{enumerate}
\subsection{Geometric structure under blow-up}
In order to get step one we give the following definition:
\begin{defin}
Let $M$ be a two dimensional complex manifold and $S_1,\cdots,S_n
\subset M$ given curves. We say that a point $p$  is a
\textbf{double intersection point} if $p$ belongs to exactly  two distinct curves among $S_1,\cdots,S_n$. If instead $p$ belongs to
 exactly three of them it is called  a \textbf{triple intersection point}.
\end{defin}
\begin{osserv}\label{osserv:o1}
In the study of curve desingularization the set of curves we find is composed by the strict transform of
the curve $S$ and the several exceptional divisors obtained by succesive blow-ups.
Because of the structure of the blow-up process we can only have double and triple intersection points
(see \cite{LAU}).
A triple intersection point belongs to the strict transform of $S$ and to two exceptional divisors.
To distinguish these two $\C\mathbb{P}^1$ we will call  \emph{ old exceptional divisor} the strict
 transform of a given exceptional divisor. Instead we will call \emph{new exceptional divisor} the
exceptional divisor produced by the last blow-up.
\end{osserv}
Now we can describe the geometric evolution  under blow-up. Note
that the only intersection point that can be triple is the one
made up by the strict transform of $S$. We will prove the
following behaviour.
\begin{propos}\label{propos:ordine}
Let $S$ be a singular curve and  let $p$ be a singularity of $S$.
The resolution process of $S$ in $p$ is related to the behavior of the multiplicity of $S$ in $p$ in the following way:
\begin{itemize}
\item
If we blow-up a singularity and the multiplicity does not reduce we have two cases:
\begin{enumerate}
\item
if we are in a double intersection point at the next blow-up we find another double intersection,
\item
if  we are in a triple intersection point at the next blow-up we can
find either a double intersection  or a triple intersection point. More precisely we
find a double intersection point if the tangent cone to the curve does
not coincide with any exceptional divisor,  while  we find a
triple intersection point if the tangent cone coincides  with one of the
two exceptional divisors and the new triple intersection point belongs to the strict transform of the
old exceptional divisor.
\end{enumerate}
\item
If we blow-up a singularity and the multiplicity reduces  we have two cases:
\begin{enumerate}
\item
if we are in a double intersection point  at the next blow-up we find a triple intersection point,
\item
if we are in a triple intersection point at the next blow-up we find a triple intersection point
that belongs to the strict transform of the new  exceptional divisor.
\end{enumerate}
\end{itemize}
\end{propos}
\begin{osserv}
In the previous Proposition we have used inproperly the expression
``the tangent cone coincides  with one of the two exceptional
divisors'' to mean that the tangent cone of $S$ in $p$ coincides
with the tangent space of $D$ in $p$.
\end{osserv}
In order to get Proposition \ref{propos:ordine} we will prove some elementary Lemmas.
\begin{lemma}\label{lemma:giallo}
Let $M$ be a two dimensional complex manifold, $S$ an analytic
irreducible curve on $M$ and $p\in S$ a singularity of $S$.
Blow-up $M$ in $p$ and let $\hat{S} $ be  the strict transform of
$S$, $D$ the exceptional divisor and $\hat{p}:=\hat{S} \cap D$.
The multiplicity of $\hat{S}$ in $\hat{p}$ is strictly smaller
than the multiplicity of $S$ in $p$ if and only if $D$ coincides
with the tangent cone of $\hat{S}$ in $\hat{p}$.
\end{lemma}
\begin{proof}
We can assume that $ p=(0,0) $ and $ S = \{ l (x,y)=0\} $ with $ l
(x,y)= y^m + l_{m+1}(x,y) + \cdots .$ Blow-up in $p$ and using the
chart such that the projection becomes  $\pi(u,v)=(u,uv)$ we have:
$ \hat{S}=\{\hat{l}(u,v)=0\} $,  with $ l(u,v)=v^m + u
l_{m+1}(1,v)+ \cdots = v^m + u q_{k-1} + \cdots $ and $ D =
\{u=0\} . $ The multiplicity of  $ \hat{S} $ in $ (0,0) $ is
strictly less then  $ m $ if and only if  $ k<m $ and then if and
only if the tangent cone  is  $ \{ u q_{k-1}(u,v)=0\} $ and so if
and only if $ D $ is included in the tangent cone. Because $S$ is
irreducible this can happen if and only if $
q_{k-1}(u,v)=u^{k-1}$, i.e. if and only if $ D $ is the tangent
cone.
\end{proof}

\begin{lemma}\label{lemma:viola}
Let $M$ be a two dimensional complex manifold, $S$ an analytic
irreducible curve on $M$ and $p\in S$ a singularity of $S$.
Blow-up $M$ in $p$ and let $\hat{S} $ be the strict transform of
$S$, $D$ the exceptional divisor and $\hat{p}:=\hat{S} \cap D$.
The exceptional divisor $D$ is  the tangent cone of $\hat{S}$ in
$\hat{p}$ if and only if blowing-up in $\hat{p}$ we get a triple
intersection point.
\end{lemma}
\begin{proof}
Let $\hat{D}$ be the strict tranform of $D$ and $D_1$ the new
exceptional divisor. Now $\hat{D}$ intersects $D_1$ in the point
corresponding to the tangent of $D$ in $p$, so $\hat{D}\cap
\hat{\hat{S}} \neq \emptyset$ if and only if $D$ and $\hat{S}$
have the same tengent in $p$. So we get a triple intersection
point if and only if the tangent cone of $\hat{S}$ coincides with
$D$.
\end{proof}
Using the previous two Lemmas we obtain the following:
\begin{lemma}
Let $ S \subset M $ be an analytic irreducible curve of
multiplicity  $ m $ in the singular point $p$. Suppose that after
a finite number of blows-up the strict transform of $S$, $
\tilde{S} $, intersects the exceptional divisor in a point
$\tilde{p}$ and indicate with $D$ the irreducible component of the
exceptional divisor conteining $\tilde{p}$, i.e. $\tilde{p}$ is a
double intersection point. Blow-up in $\tilde{p}$ and let $ D_1 $
be the new exceptional divisor and  $\hat{S} $ the strict
transform of $ \tilde{S} $. If the multiplicity of $ \hat{S} $ in
$ \hat{p}:=D_1\cap \hat{S}  $ is equal to the multiplicity of
$\tilde{S} $ in $\tilde{p}$ then at the following blow-up we find
again a double intersection point.
\end{lemma}
By Lemma \ref{lemma:viola} we also get:
\begin{lemma}\label{lemma:verde}
Let $ S \subset M $ be an analytic irreducible curve of
multiplicity  $ m $ in the singular point $p$. Suppose that after
a finite number of blows-up  we have a triple intersection point.
At the following blow-up we have two cases:
\begin{enumerate}
\item  if the tangent cone in the singularity contains one of the two exceptional divisors then at the next blow-up we find agin a triple
intersection point,
\item if the tangent cone in the singularity does not contain any of the two exceptional divisors then at the next blow-up we find a
 double intersection point.
\end{enumerate}
\end{lemma}
\begin{osserv}\label{osserv:blu}
We observe that the demonstrative method used in Lemma
\ref{lemma:verde} does not give informations on which of the
exceptional divisors goes to create the new triple intersection.
To get this information we need some more calculations. Let
$\hat{S}$ the strict transform of $S$ after some blow-ups and
suppose to have a triple intersection point.
 We can assume that $p=(0,0) $ and $ \hat{S}=\{\hat{l}(u,v)=0\} $ with  $\hat{l}(u,v)=
 v^{m} + u^{k_1}[q_{k_2-k_1}(u,v)+\cdots] ,$ and  $ D_1=\{v=0\} $ , $ \hat{D}=\{u=0\}$
 where $D_1$ is the new exceptional divisor and $\hat{D}$ is the
 old one (according to Remark \ref{osserv:o1}).
Let esaminate the various cases:
\begin{enumerate}
\item  If $ m>k_2 $ then the tangent cone is $
\{u^{k_1}q_{k_2-k_1}(u,v)=0\} $ and by the irreducibility of $S$
is $ \{cu^{k_2}=0\} $ with $ c\neq 0 $ and so it contains an
exceptional divisor, $ \hat{D} .$ Blow-up again  $ (0,0) $ and
using the chart by which the projection is $ \pi(x,y)=(xy,y)$ we
have:
$$ \hat{l}(xy,y)= y^{m} + c x^{k_2}y^{k_2} + x^{k_1 }y^{k_2 +1}[ q_{k_2-k_1+1}+ \cdots ] $$
and because  $ m>k_2 $
$$ \hat{\hat{l}}(x,y)=y^{m_1-k_2}+ c x^{k_2} + x^{k_1} y[q_{k_2-k_1+1}(x,1)+\cdots ]$$
with $ D_2=\{y=0\}$ e $ \hat{\hat{D}}=\{x=0\} .$ So $ (0,0)$ is a
triple intersection point made up by $D_2$, $\hat{\hat{S}}$,
$\hat{\hat{D}}$. If instead we use the other chart we find only a
double intersection points. \item If $ m_1<k_2 $ we proceed in the
same way obtaining a triple intersection point made by $D_2$,
$\hat{D_1}$ and $\hat{\hat{S}}$. \item If  $ m_1=k_2 $  the
tangent cone is given by $ \{v^{m_1}+u^{k_1}q_{k_2-k_1}(u,v)=0\} $
and by the irreducibility of the curve it is $ \{(v+c u)^{m_1}=0\}
$ with $ c\neq 0 $ and it does not contain any exceptional
divisor. So by Lemma \ref{lemma:verde} at the next blow-up we find
only double intersection points.
\end{enumerate}
\end{osserv}

\subsection{C.S.S. index under blow-up}
Now we can proceed in order to get step two studying the behaviour
of the index in a general resolution process via blow-up. The
upshot is to prove that in the resolution process we necessarily
find  a $(C_1)$ or $(C_2)$ point ,i.e., $p$ is an appropriate
singularity and then Theorem \ref{teor:miomio} holds.

The intent is to analyze the C.S.S. index in all possible geometric evolutions (see Proposition~\ref{propos:ordine}).
\begin{osserv}
In the analysis we will omit the case in which at some blow-up we
find a dicritical point (see Definition 3.2 in \cite{BR}). In fact
in this case the goal is obtained by Proposition 7.8  \cite{BR}
and by the proper mapping theorem  \cite{GR}.
\end{osserv}
We will consider resolution processes only at a combinatoric level in a sense that will be specified later.

Thanks to Proposition \ref{propos:ordine} the structure of a
resolution process of a singular point $p$ is completely described
by the behaviour of the multiplicity of the strict transform at the
intersection with the exceptional divisor. We can then  consider a
sequence of blow-ups only as a sequence of positive number
(representing the evolution of the multiplicity) and forgetting
any type of geometric obstruction.
\begin{defin}
A \textbf{process} is  an ordinate  list of the form:
$$ P=\{(k, m), (\alpha_1, m_1), \cdots, (\alpha_n ,m_n)\} $$
where $k,\alpha_i, m_i \in \mathbb{N}$ and $ m> m_1\geq \cdots
\geq m_n$. We associate to P, from a purely formal  point of view,
a  blow-up sequence for a curve $S$ where the blows-up are made at the beginning at 
the point $p$  and then at the intersection point of the strict transform of the curve and the  exceptional divisor.
The blow-up sequence satisfies the following rules:
\begin{itemize}
\item[-]from the first to the $k-th$ blow-up we find only double intersection points and the  curve multiplicity 
is constantly equal to  $m$,
\item[-]
from the  $(k+1)-th$ to the $(k+\alpha_1)-th$ blow-up  we find a triple intersection point
and the multiplicity of the strict transform of $S$ is constantly equal to $ m_1<m $,
 \item[]  \hspace{5cm} \vdots
\hspace{3cm}  \vdots 
\item[-]from the  $(k+\alpha_1+\cdots +\alpha_{n-1}+1)-th$ to the $(k+\alpha_1 + \cdots \alpha_n)-th$ blow-up
we find a triple intersection  point and the  multiplicity is constantly equal to 
$m_n \leq m_{n-1}$.
\end{itemize}
\end{defin}
\begin{osserv}
At the end of $P$ the curve $S$ is not desingularized, in fact we
have triple points and this type of point are not admitted in the
desingularized curve.
\end{osserv}
Now, according to Proposition \ref{propos:ordine} we start to
analyze all the possible cases. For notations we refer to
\cite{BR} and \cite{BRU}.
\subsection{Case of double intersection}
It corresponds to a process $ P=\{(k, m)\}$, i.e. we start with
multiplicity $m$ and we remain with this multiplicity for $k$
blows-up finding only double points. If we do not find $(C_1)$ or
$(C_2)$ points in the total transform then (arguing as in
Proposition 7.8(2) of \cite{BR}) at the $k-$th blow-up the indices
are of type:
\begin{equation}\label{eq:doppia1}
\begin{aligned}  \id(\tilde{\F},D,q)&\in \Q_{\leq-\frac{1}{k}}\\
                   \id(\tilde{\F},\hat{S},q)&\not\in \Q_{\geq - k m^2}.
\end{aligned}
\end{equation}
where  $q:= \hat{S}\cap D $.
\subsection{Case of triple intersection}
We consider now a slightly more complicated process,\\
$P=\{(k,m),(1,m_1),\cdots,(1,m_n)\}$. Let us suppose not to find
$(C_1)$ or $(C_2)$ points during $P$.

We indicate at the last blow-up with  $S$ the strict transform of
the curve, \Fl the saturated foliation, $D_1$, $D_2$ the two
exceptional divisors that intersect, with $S$, in the last triple
intersection point $q$.

\begin{propos}\label{propos:ciao}
In this situation at the last blow-up of $P$, if we have not found $(C_1)$ or $(C_2)$ points, we can find
$ x,y\in\mathbb{N} $ and  $ a,b \in\mathbb{N}\cup\{0\}$ such that the indices are:
\begin{equation} \label{eq:i2}
\begin{aligned}
&\id(\F,S,q)\not\in \Q_{\geq-km^2- m_1^2 - \cdots -  m_n^2}\\
&\id(\F,D_1,q)\in \Q_{\leq -\frac{x}{y}}\\
&\id(\F,D_2,q)\in \Q_{\leq -\frac{y k + a}{x k + b}}.
\end{aligned}
\end{equation}
\end{propos}
\begin{proof}
At the k-th blow-up the indices are of type \eqref{eq:doppia1}. Let  blow-up again.
As P describes we have a multiplicity decrease and we find a triple point on the total transform.
Then if some point of the new exceptional divisor $D_1$ is of type $(C_1)$, $p$ is an appropriate singularity and we have the assertion.
Otherwise  $\id(\F,D_1,p)\in\Q_{\geq 0}$ $ \forall p\in D_1\setminus \{q\}$ and then by Index Theorem:
\begin{equation*}
\id(\F,D_1,q) \in \Q_{\leq-1}.
\end{equation*}
Then by Proposition \ref{propos:blow-up} and observing that $D$ has multiplicity one:
\begin{equation*}
\begin{aligned}
&\id(\F,\hat{S},q) \not \in \Q_{\geq -k m^2 - m_1^2}\\
&\id(\F,\hat{D_2},q)\in \Q_{\leq-\frac{k+1}{k}}.
\end{aligned}
\end{equation*}
Proceeding by induction on $n$ we can assume the assertion true for $n$ and we prove it for $n+1$.
We have to analyze separately  two different cases that can occur blowing-up:
\begin{enumerate}
                         \item the new triple point is made by  $ \{\hat{S},\hat{D}_2,D\} ;$
                         \item the new triple is made by  $ \{\hat{S},\hat{D}_1,D\},$
\end{enumerate}
where $D$ is the new exceptional divisor and $D_1$ and $D_2$ are
the ones of the $n$ blow-up whose indices satisfy \eqref{eq:i2} by
inductive hypothesis. We consider only the case $(1)$ because the
other is similar. By Proposition \ref{propos:blow-up} the indices
are of type:
\begin{equation*}
\begin{aligned}
&\id(\F,\hat{S},q_1)=\id(\F,S,p)-m_1^2- \cdots -m_n^2 - m_{n+1}^2\\
&\id(\F,\hat{D}_2,q_1)\in \Q_{\leq -\frac{(x+y)k+(a+b)}{xk+b}}\\
&\id(\F,\hat{D}_1,q_0)\in \Q_{\leq -\frac{x+y}{y}},
\end{aligned}
\end{equation*}
where $q_1$ is the new triple point and $q_0:=\hat{D_1} \cap D$.
If there are not  $ (C_1) $ points on  $ D\setminus\{q_0,q_1\} $ then by Index Theorem $ q_0 $ is a $ (C_2) $ point or
$ \id(\F,D,q_1)\in \Q_{\leq -\frac{x}{x+y}} .$
In the last case the indices satisfy:
\begin{equation}
\begin{aligned}
&\id(\F,\hat{S},q_1)=\id(\F,S,p)-m_1^2- \cdots -m_n^2 - m_{n+1}^2\\
&\id(\F,\hat{D}_2,q_1)\in \Q_{\leq -\frac{(x+y)k+(a+b)}{xk+b}}\\
&\id(\F,D,q_1)\in \Q_{\leq -\frac{x}{x+y}}.
\end{aligned}
\end{equation}
and then the assertion follows putting $ y'= x+y, x'=x, a'=a+b, b'=b.$
\end{proof}

\begin{osserv}\label{osserv:proc}
A general process can always be wrietten in the form $P=\{(k,m),(\alpha_1,m_1),\cdots
,(\alpha_n,m_n)\}$ with $m_i \neq m_j$ if $i\neq j$.
The coefficients $(x,y,a,b)$ that occour in $P$, by Proposition \ref{propos:ciao} depend  only on  the  $ \alpha_i $ and 
to the order in which they appear but  not to the multiplicities $m_i$ and
the coefficient $k$.
\end{osserv}
We propose now some simple properties of the index under a process that will be usefull later:
\begin{lemma}\label{lemma:l0}
In \eqref{eq:i2} it follows that  $ x a - y b = 1. $
\end{lemma}
\begin{proof}
We proceed by induction on the number of blows-up and argue as in the proof of Proposition \ref{propos:ciao}
\end{proof}
With the same arguments we can also prove:
\begin{lemma}\label{lemma:l1}
Let consider a process $P=\{(k,m),(1,m_1),\cdots,(1,m_n)\}$ and
indicate with  $S, D_1, D_2$ the curves that create the triple
intersection point. Then if  $ (x,y,a,b) $ are the coefficients
that appear in the indices \eqref{eq:i2}  we have, according to
Remark \ref{osserv:o1}:
\begin{itemize}
\item[] if $ x > y $ then $ D_2 $  is the new exceptional divisor and   $  D_1$ is the old one,
\item[] if $ x \leq y $ then  $ D_1 $ is the new exceptional divisor and  $  D_2 $  is the old one.
\end{itemize}
\end{lemma}
Using Lemma \ref{lemma:l1} and Remark \ref{osserv:blu} we can
easily prove:
\begin{lemma}\label{lemma:somma}
If we blow-up a triple intersection point and we have a multiplicity decrease  then the coefficients
$ (x',y',a',b')$ of the indices of the new triple are such that:
\begin{itemize}
\item[] if $ x > y $ then $ x'=x, y'=x+y ,$
\item[] if $ x \leq y $ then  $ x'=x+y, y'=y .$
\end{itemize}
\end{lemma}
In the analysis of the $C.S.S.$ index in the triple intersection case the knowledge of the index is equivalent to the knowledge
 of the coefficients $(x,y,a,b)$.
According to Remark \ref{osserv:blu} the decrease or not of the multiplicity gives different coefficients.
In the next subsections we are going to investigate these cases.
To make clearer the possible evolutions of the coefficients we report below the coefficients $(x,y,a,b)$ that can appear in the first
 five blows-up in triple intersection.
We indicate in black the coefficients related to a decrease of multiplicity and in grey the others.
\begin{figure}[h] 
\begin{center}
\includegraphics[width=1   \textwidth]{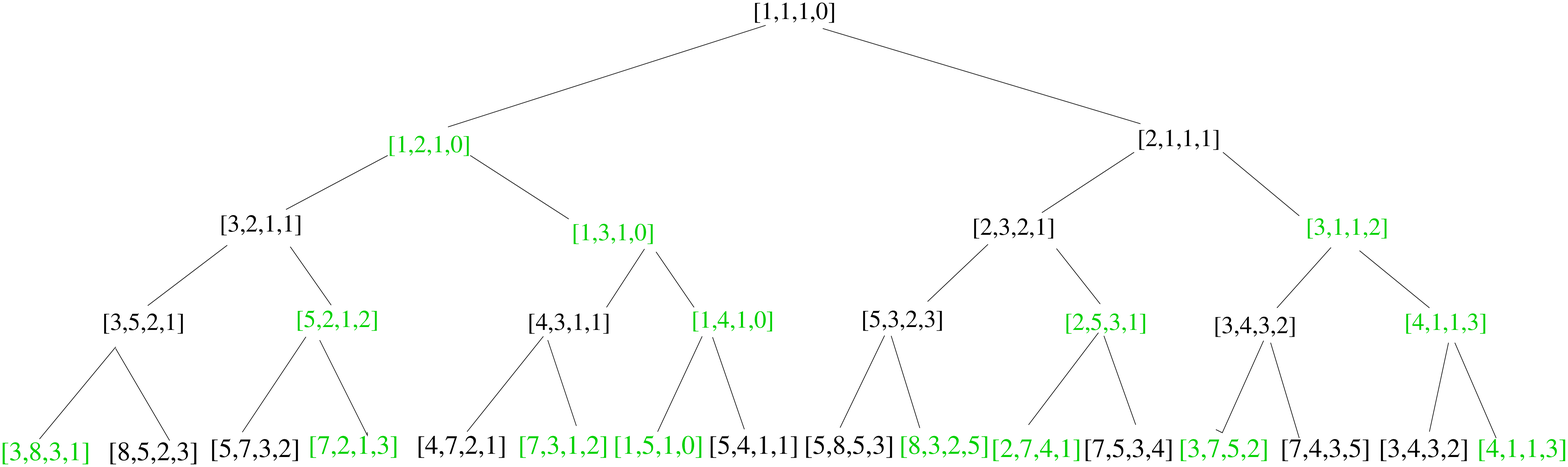}
\end{center}
\end{figure}
\subsection{Transition from triple intersection with multiplicity lowering to triple with
constant multiplicity} We consider a process of type
$P=\{(k,m),(\alpha_1, m_1),\cdots,
(\alpha_{n-1},m_{n-1}),(\alpha_n, m_n)\}$ with $m_i\neq m_j $ if
$i\neq j.$ We want to relate the coefficients of the last blow-up
with the ones obtained at the first lowering of multiplicity
$m_{n-1}\rightarrow m_n$, i.e.,  we want to relate the last
indices of the process $\{(k,m),(\alpha_1, m_1),\cdots,
(\alpha_{n-1},m_{n-1}),(1, m_n)\}$ to the last ones of $P$.
\begin{propos}\label{section:TT}
Suppose to have indices of type:
\begin{equation}\label{eq:in}
\begin{aligned}
&\id(\F,S,p)\not \in \Q_{\geq -k m^2- \alpha_1 m_1^2- \cdots -\alpha_{n-1}m_{n-1}^2-m_n^2}\\
&\id(\F,D_1,p)\in \Q_{\leq{-\frac{x}{y}}}\\
&\id(\F,D_2,p)\in \Q_{\leq-\frac{y k + a}{x k + b}}
\end{aligned}
\end{equation}
with  $ m_i\neq m_j $ if $i\neq j $, i.e.,  $ n $  is the number of multiplicity lowerings.
The indices at the end of the process $P$ are of type:
\begin{itemize}
\item[]
if $ x>y $
\begin{equation}\label{eq:in1}
\begin{aligned}
&\id(\F,S,p)\not \in \Q_{\geq -k m^2- \alpha_1 m_1^2- \cdots -\alpha_{n-1}m_{n-1}^2-\alpha_{n}m_n^2}\\
&\id(\F,D_1,p)\in \Q_{\leq{-\frac{x+(\alpha_n-1)y}{y}}}\\
&\id(\F,D_2,p)\in \Q_{\leq-\frac{y k + a}{(x+(\alpha_n-1)y) k +((\alpha_n-1)a+ b)}}
\end{aligned}
\end{equation}
\item[]
if $ x\leq y $
\begin{equation}\label{eq:in2}
\begin{aligned}
&\id(\F,S,p)\not \in \Q_{\geq -k m^2- \alpha_1 m_1^2- \cdots -\alpha_{n-1}m_{n-1}^2-\alpha_{n}m_n^2}\\
&\id(\F,D_1,p)\in \Q_{\leq{-\frac{((\alpha_n-1)x+y)k + (a + (\alpha_n-1)b)}{x k+y}}}\\
&\id(\F,D_2,p)\in \Q_{\leq-\frac{x}{(\alpha_n-1)x+y}}.
\end{aligned}
\end{equation}
\end{itemize}
\end{propos}
\begin{proof}
By Proposition \ref{propos:ordine} blowing-up with constant  multiplicity we know that the new triple point is made up by the curve, the
new exceptional divisor and the strict transform of the old one (see Remark \ref{osserv:o1}). 
We have to analyze separately the case in which the old exceptional divisor is $D_1$ or $D_2.$ This distinction can
be made in terms of $x>y$ or $x\leq y$ thanks to Lemma \ref{lemma:l1}.
Suppose, for instance, $ x>y $ in the indices  \eqref{eq:in},  then we conclude that the old exceptional divisor is $D_1$.
Now  blowing-up again and using Proposition \ref{propos:blow-up}, the Index theorem and the assumption of non existence of
 $(C_1)$ or $(C_2)$  points we can prove the result for $\alpha_n=1,2$.
Then proceeeding by induction and repeating the same argument  for the case $x\geq y$ we have the assertion.
\end{proof}

\subsection{Transition from triple to double intersection}\label{subsection:td}
Suppose that, after $ k $ blows-up in double intersection and a
finite number of blows-up in triple intersection,  we return to
double intersection. Let consider the  generic indices of the
triple  \eqref{eq:i2} and we write the index along $S$ in  the
form:
\begin{equation*}
\id(\F,S,p)\not \in \Q_{\geq -k m^2-\alpha_1 m_1^2-\cdots -\alpha_n m_n^2},
\end{equation*}
with $m_i\neq m_j$ if $i\neq j$. Using Lemma \ref{lemma:l0} we obtain
that the indices in the double point we find  are:
\begin{equation}\label{eq:doppia}
\begin{aligned}
&\id(\F,S,q)\not \in \Q_{\geq-k m^2-\alpha_1 m_1^2-\cdots -\alpha_n m_n^2-m_n^2}\\
&\id(\F,D,q)\in \Q_{\leq-\frac{1}{(x+y)^2 k + (x+y)(a+b)}}.
\end{aligned}
\end{equation}

\subsection{Estimate of the term  $ k m^2 $}\label{subsection:fine1}
We estimate the term $-k m^2-\alpha_1 m_1^2-\cdots -\alpha_n m_n^2-m_n^2 $ showing that, if the curve is resolved, then $q$ is
a point of type $(C_2)$; otherwise we obtain indices of the form \eqref{eq:doppia1} and so we can
utilize again the results found in the previous sections in order to get  desingularization.
In this subsection we estimate the term $k m^2 $.
\begin{propos}
If we indicate with  $ (x_i^j,y_i^j,a_i^j,b_i^j) $ the
coefficients that appear in the indices of the triple intersection point at the
$j-$th blow-up with multiplicity $ m_i $ then, if $ n \geq 2 $:
\begin{eqnarray}\label{eq:st1}
m&=& x_{n-1}^{\alpha_{n-1}} m_{n-1}+ y_{n-1}^{\alpha_{n-1}} m_n \mbox{ if } x_{n-1}^{\alpha_{n-1}}\geq
 y_{n-1}^{\alpha_{n-1}},\\\label{eq:st2}
m&=& y_{n-1}^{\alpha_{n-1}} m_{n-1}+ x_{n-1}^{\alpha_{n-1}} m_n \mbox{ if } y_{n-1}^{\alpha_{n-1}}\geq
 x_{n-1}^{\alpha_{n-1}}.
\end{eqnarray}
\end{propos}
\begin{proof}
We proceed by induction on the number $ n $ of changes of
multiplicity. For $n=2$ the indices are of the form:
\begin{equation*}
\begin{aligned}
&\id(\F,S,p)\not\in \Q_{\geq- k m^2-\alpha_1 m_1^2-\alpha_2m_2^2},\\
&\id(\F,D_1,p)\in\Q_{\leq -\frac{\alpha_1\alpha_2+1}{\alpha_1}},\\
&\id(\F,D_2,p)\in\Q_{\leq-\frac{\alpha_1 k +
1}{(\alpha_1\alpha_2+1)k + \alpha_2}}.
\end{aligned}
\end{equation*}
The indices we find at the $\alpha_1-$th blow-up with multiplicity $ m_1 $ are:
\begin{equation}\label{eq:n}
\begin{aligned}
&\id(\F,S,p)\not\in \Q_{\geq - k m^2-\alpha_1 m_1^2},\\
&\id(\F,D_2,p)\in\Q_{\leq -\frac{\alpha_1k+1}{k}},\\
&\id(\F,D_1,p)\in\Q_{\leq-\frac{1}{\alpha_1}}.
\end{aligned}
\end{equation}
Because we make  $ \alpha_1 $ blows-up with multiplicity $ m_1 $ and because the curve is irreducible by Enriques-Chisini theorem
(\cite{Bri} pag. 516) we have:
$$ m_2=m-\alpha_1 m_1 $$
and then the assertion. We prove the inductive step.
The index along $S$ is:
\begin{equation*}
\id(\F,S,p)\not\in\Q_{\geq -k m^2-\alpha_1 m_1^2-\cdots
-\alpha_{n-1}m_{n-1}^2-\alpha_n m_n^2-\alpha_{n+1}m_{n+1}^2}
\end{equation*}
We consider  the case  $ x_{n-1}^{\alpha_{n-1}}\geq y_{n-1}^{\alpha_{n-1}}$  ( the other is similar ).
Because we make $ \alpha_n $ blows-up with multiplicity  $ m_n $ we have:
$$ m_{n+1}=m_{n-1}-\alpha_n m_n \mbox{ and then  } m_{n-1}= \alpha_n m_n + m_{n+1}. $$
By inductive hypothesis and the above relation we find an expression of $m$
in terms of $m_n$ and $m_{n+1}$.
Now we have to prove that this expression is the one of the statement.
Using Lemma \ref{lemma:somma} we have that $ x_n^1 \leq y_n^1 $
and for Proposition  \ref{section:TT} the indices at the $ \alpha_n-$th blow-up with
multiplicity  $ m_n $ are:
\begin{equation*}
\begin{aligned}
&\id(\F,S,p)\not\in \Q_{\geq -k m^2-\alpha_1 m_1^2- \cdots -\alpha_{n-1}m_{n-1}^2-\alpha_{n}m_n^2},\\
&\id(\F,D_2,p)\in \Q_{\leq - \frac{(y_n^1+(\alpha_n-1)x_n^1)k +
(a_n^1+(\alpha_n-1)b_n^1)}
{x_n^1 k + b_n^1}},\\
&\id(\F,D_1,p)\in\Q_{\leq -\frac{x_n^1}{(\alpha_n-1)x_n^1+y_n^1}}.
\end{aligned}
\end{equation*}
Clearly $ x_n^{\alpha_n}\leq y_n^{\alpha_n}$ and so computing the expression
 $ y_n^{\alpha_n}m_n + x_n^{\alpha_n}m_{n+1},$ using the above form of the coefficients
and Lemma \ref{lemma:somma} we get the assertion.
\end{proof}
\begin{propos}\label{propos:c1}
When in the resolution process we return in double intersection the indices:
\begin{equation}\label{eq:i3}
\begin{aligned}
&\id(\F,S,p)\not \in \Q_{\geq -k m^2-\alpha_1m_1^2-\cdots -\alpha_nm_n^2-m_n^2},\\
&\id(\F,D,q_0)\in \Q_{\leq-\frac{1}{(x_n^{\alpha_n}+y_n^{\alpha_n})^2 k +
(x_n^{\alpha_n}+y_n^{\alpha_n})(a_n^{\alpha_n}+b_n^{\alpha_n})}},
\end{aligned}
\end{equation}
satisfy
$$ m \geq (x_n^{\alpha_n}+y_n^{\alpha_n})m_n .$$
\end{propos}
\begin{proof}
It derives directly from the previous proposition and from Lemma
\ref{lemma:somma}.
\end{proof}

\subsection{Estimate of the terms  $k m^2+ \alpha_1 m_1^2\cdots + \alpha_n m_n^2+ m_n^2 $}\label{subsection:fine2}
\begin{propos}\label{propos:restanti}
The indices at the return in double intersection \eqref{eq:i3}, with $ n\geq 2$, satisfy:
$$ \alpha_1 m_1^2+\cdots \alpha_n m_n^2+ m_n^2 \geq
 (x_n^{\alpha_n}+y_n^{\alpha_n})(a_n^{\alpha_n}+b_n^{\alpha_n})m_n^2 .$$
\end{propos}
Before proving this statement we consider the following one:
\begin{propos}\label{propos:processi}
Let  $ P=\{(k,m), (\alpha_1, m_1), \cdots, (\alpha_n, m_n)\} $ be
a process and let indicate with
 $ (x,y,a,b) $ the coefficients of the indices that appear at the last blow-up described by $P$.
We associate to $P$ the process  $ \bar{P}=\{(k, m),(\alpha_2,
m_2),\cdots,(\alpha_n, m_n)\} $ and we indicate with  $
(\bar{x},\bar{y},\bar{a},\bar{b})$ the coefficients of the indices
that appear at the last blow-up described by $\bar{P}$. Then:
\begin{equation*}
\begin{aligned}
b&=\bar{y} \hspace{3cm} x=\alpha_1\bar{y}+\bar{a} \\
a&=\bar{x} \hspace{3cm} y=\alpha_1\bar{x}+\bar{b}
\end{aligned}
\end{equation*}
\end{propos}
\begin{proof}
We proceed by induction on the number $n$ of multiplicities
decreases. By a direct calculation the Proposition is true for
$n=2$.
Suppose the assertion true for $ n $ and let prove it
for $n+1$.\\
 Let consider the two processes $ P'=\{(k,m),(\alpha_1,m_1),\cdots,(\alpha_n, m_n),(\alpha_{n+1}, m_{n+1})\} $
and \\ $ \bar{P}'=\{(k,m),(\alpha_2, m_2), \cdots,
(\alpha_n,m_n),(\alpha_{n+1},m_{n+1})\} $ with respectively end
coefficients
$ (x',y',a',b') $ and $ (\bar{x}',\bar{y}',\bar{a}',\bar{b}'). $\\
Let now construct the following two processes
$ P=\{(k,m),(\alpha_1,m_1),\cdots,(\alpha_n, m_n)\},$ \\
$ \bar{P}=\{(k,m),(\alpha_2,m_2),\cdots,(\alpha_n, m_n)\} $ with
end coefficients  $ (x,y,a,b) $ and
$(\bar{x},\bar{y},\bar{a},\bar{b}) .$ Starting by coefficients  $
(x,y,a,b) $ we get  $ (x',y',a',b') $ after one blow-up with
multiplicity decrease and other $ \alpha_{n+1}-1 $ blows-up with
constant multiplicity $m_{n+1}.$ By Propositions
\ref{propos:ordine} and \ref{section:TT}:
\begin{equation*}
\begin{aligned}
(x',y',a',b')&=(x,y+\alpha_n x,a+\alpha_n b,b) \mbox{ if } x>y,\\
(x',y',a',b')&=(x+\alpha_n y,y,a,\alpha_n a+b) \mbox{ if } x \leq y.
\end{aligned}
\end{equation*}
Similarly we get:
\begin{equation*}
\begin{aligned}
(\bar{x}',\bar{y}',\bar{a}',\bar{b}')&=(\bar{x},\bar{y}+\alpha_n \bar{x},\bar{a}+\alpha_n \bar{b},\bar{b})
 \mbox{ if } \bar{x}>\bar{y},\\
(\bar{x}',\bar{y}',\bar{a}',\bar{b}')&=(\bar{x}+\alpha_n \bar{y},\bar{y},\bar{a},\alpha_n \bar{a}+\bar{b})
 \mbox{ if } \bar{x}\leq\bar{y}.
\end{aligned}
\end{equation*}
The processes  $ P $ e $ P' $ differs only on one multiplicity
decrease.  Propositions \ref{propos:ordine} and \ref{section:TT}
say that $ x $ and $ y $ relations invert only when a multiplicity
decrease occurs. Then we can conclude that $ \bar{x}>\bar{y} $ if
and only if $ x \leq y .$ If, for instance, $ x> y ,$ by inductive
hypothesis:
\begin{equation*}
\begin{aligned}
b'&=b=\bar{y}=\bar{y}',\\
a'&=a+\alpha_n b= \bar{x}+\alpha_n \bar{y}=\bar{x}',\\
x'&=x=\alpha_1 \bar{y}+\bar{a}=\alpha_1\bar{y}'+\bar{a}',\\
y'&=y+\alpha_n x=\alpha_1\bar{x}+\bar{b}+\alpha_n\bar{y}+\alpha_n\bar{a}=\alpha_1(\bar{x}+\alpha_n \bar{y})+
(\bar{b}+\alpha_n \bar{a})=\alpha_1\bar{x}'+\bar{b}'.
\end{aligned}
\end{equation*}
and then the assertion.
\end{proof}

Now we can prove Proposition  \ref{propos:restanti}.
\begin{proof}
Let proceed by induction on the number of changes of multiplicity.
If  $ n=2 $  the structure of the indices can be easily computed
to obtain the assertion. Let prove the inductive step. Let
$P=\{(k,m),(\alpha_1,m_1),\cdots,(\alpha_n,m_n),(\alpha_{n+1},m_{n+1})\}$
be a generic process.
 Thanks to the inductive step applied on the process
$ \bar{P}=\{(k,m),(\alpha_2,m_2),\cdots,(\alpha_{n+1},m_{n+1})\} $
we have:
\begin{equation*}
\alpha_1m_1^2+\cdots\alpha_n m_n^2+ \alpha_{n+1} m_{n+1}^2 \geq \alpha_1 m_1^2 +
 (\bar{x}+\bar{y})(\bar{a}+\bar{b})m_{n+1}^2.
\end{equation*}
In order to estimate $ \alpha_1 m_1^2 $ we consider the process $
P'=\{(\alpha_1,m_1),(\alpha_2,m_2),\cdots,(\alpha_{n+1},m_{n+1})\}
$ and thanks to  Remark \ref{osserv:proc} and Proposition
\ref{propos:c1}  we have:
$$ m_1^2 \geq (\bar{x}+\bar{y})^2 m_{n+1}^2. $$
Then:
\begin{equation*}
\begin{aligned}
\alpha_1m_1^2+\cdots\alpha_n m_n^2+ \alpha_{n+1} m_{n+1}^2 & \geq \alpha_1 (\bar{x}+\bar{y})^2m_{n+1}^2+
(\bar{x}+\bar{y})(\bar{a}+\bar{b})m_{n+1}^2\\
& =(\bar{x}+\bar{y})(\alpha_1 \bar{x}+\alpha_1 \bar{y} +
\bar{a}+\bar{b})m_{n+1}^2.
\end{aligned}
\end{equation*}
We conclude thanks to  Proposition \ref{propos:processi}.
\end{proof}
\begin{osserv}
The estimate of $km^2$ and of the remaining terms are valid only
if $n\geq 2$. The case $n=1$ can be easily proved using
equation \eqref{eq:n}, Section \ref{subsection:td} and observing
that because of the $\alpha_1 + 1 $ blows-up $m \geq (\alpha_1 +1)m_1$.
\end{osserv}

\subsection{Proof of the Theorem}

All the previous separate particular cases can now be glued together to get Theorem \ref{teor:miomio}.
We have observed that in the resolution process we can have only double or triple intersection points and so we studied the index in
these cases.

The triple point case presents two  different subcases, linked to
the multiplicity of the curve: it can decrease or not. This
information is extremely useful for the study of the index
evolution because it identifies  the right exceptional divisor
that will occur in the next triple point. Now we observe that if
at the end of a process $P=\{(k,m),(\alpha_1, m_1),\cdots,
(\alpha_n, m_n)\}$ we find a double point and the curve is
desingularized  we are in the geometric conditions of a $(C_2)$
point. The indices are the ones given by equation
\eqref{eq:doppia} and by Propositions \ref{propos:c1} and
\ref{propos:restanti} we can estimate them in such a way they
became:
\begin{equation*}
\begin{aligned}
&\id(\Fl,S,p)\not \in \Q_{\geq -[(x+y)k^2+(x+y)(a+b)]}\\
&\id(\Fl,D,p)\in \Q_{\leq \frac{1}{(x+y)^2k + (x+y)(a+b)}}
\end{aligned}
\end{equation*}
and so $p$ is a $(C_2)$ point. Otherwise we are not in the right
geometric conditions, i.e. the resolution is not just ended, but
the same propositions gives indices:
\begin{equation*}
\begin{aligned}
&\id(\Fl,S,p)\not \in \Q_{\geq -[(x+y)k^2+(x+y)(a+b)]m_n^2}\\
&\id(\Fl,D,p)\in \Q_{\leq \frac{1}{(x+y)^2k + (x+y)(a+b)}}
\end{aligned}
\end{equation*}
and so we have indices exactly of the form of the ones associated to a process $P=\{h m^2\}$
and then we can apply all the previous argument to the new process which  is starting.
Such process terminates after a finite number of blows-up by theorem of resolution of singularities \cite{LAU}.

\section{Applications}
The demonstrative method used to get Theorem \ref{teor:miomio} allows us to generalize
to the case in which we start with  more than one  separatrix:
\begin{propos}\label{propos:mio2}
Let \M be  a two dimensional complex manifold,  \Fl an holomorphic foliation on  \M  and
$ S_0,S_1,\cdots,S_n $ separatrices  of  \F  passing through  a point  $ p\in M $.
Let assume that  $ S_1,\cdots,S_n $ are  non singular  and transverse each other and to $ S_0$.
If,  besides,  the indices are of the sequent form:
\begin{equation*}
\begin{aligned}
&\id(\F,S_0,p)\not\in \Q_{\geq -m^2}\\
&\id(\F,S_i,p)\in\Q_{\leq-(2n-1)} \quad \forall i\geq 1,
\end{aligned}
\end{equation*}
then another separatrix through  $p$ exists.
\end{propos}
\begin{proof}
We prove  that  $p$ is an appropriate singularity. We observe
that after the first blow-up, if we have not finished, we have the
same indices found in the study made to prove
Theorem \ref{teor:miomio} and so we conclude with the same
argument.
\end{proof}

We show briefly that Theorem \ref{teor:miomio} includes as
particular cases the classical results in discrete and continuous
dynamics.

\begin{corol}[\cite{CS}]\label{corol:mio1}
Let  \M  be a two dimensional complex manifold,  \Fl an
holomorphic foliation on \M and  $ p\in M $ a singularity of \F.
Then a separatrix of  \Fl for $p$  exists.
\end{corol}
\begin{proof}
We blow-up $M$ in $p$. If the exceptional divisor is not  a separatrix
for the saturated foliation, Proposition $1$ in \cite{BRU}$($pag.15$)$
concludes. Otherwise  using the index theorem $($see \cite{SW}$)$
and remembering that $D\cdot D =-1$, we obtain the existence of a
singularity  $ \tilde{p} $  of  the saturated foliation
$ \tilde{\F} $ such that $ \id(\tilde{\F},D,\tilde{p})\not \in
\Q^+\cup\{0\}$  and then by Theorem \ref{teor:miomio} we have the
existence of another separatrix for $ \tilde{p} $ that projects in
a separatrix for \Fl in \p.
\end{proof}
With similar arguments we also have:
\begin{corol}[\cite{SB}]
Let  $M $ be  a two dimensional complex manifold, \Fl an holomorphic foliation on  $M$.
Let  $ S\subset M $ be  a compact, globally and locally irreducible curve with  $S \cdot S<0$.
If  $S$ is a separatrix for  \Fl then  a point $ p\in S$ for which passes another separatrix exists.
\end{corol}

\begin{osserv}
Analogously to what said for Theorem \ref{teor:diff} we can
obtain, in diffeomorphisms dynamics,  a similar result to
Proposition \ref{propos:mio2} and find as particular cases results of Abate \cite{AB}
and Bracci \cite{BR}.
\end{osserv}
\begin{propos}Let \M be  a two dimensional complex manifold, $f:M
\longrightarrow M$ an holomorphic map on $M$ with $ \hbox{Fix}(f)=
S_0 \cup, S_1,\cdots,\cup S_n $ with $ S_0,\cdots,S_n $ analytic
curves passing through the same point $p \in M$. Let suppose that
$S_1,\cdots,S_n $ are non singular and transverse each other and
to $ S_0$. If, besides, the indices are of the sequent form:
\begin{equation*}
\begin{aligned}
&\id(f,S_0,p)\not\in \Q_{\geq -m^2}\\
&\id(f,S_i,p)\in\Q_{\leq-(2n-1)} \quad \forall i\geq 1,
\end{aligned}
\end{equation*}
then a parabolic curve through  $p$ exists.
\end{propos}
\begin{corol}[\cite{AB}]
Let  \M  be a two dimensional complex manifold,  $f:M
\longrightarrow M$  an holomorphic map on \M and  $ p\in M $ an
isolated singularity of $f$ such that $df_p=Id$. Then a parabolic
curve for $f$ through $p$ exists.
\end{corol}
\begin{corol}[\cite{SB}]
Let  $M $ be  a two dimensional complex manifold, $ f:M
\longrightarrow M$  an holomorphic map on  $M$. Let  $ S\subset M
$ be a compact, globally and locally irreducible curve with  $S
\cdot S<0$. If $\hbox{f}=S$ and $f$ is non degenerate along $S$ is
 then a point $ p\in S$ for which passes a parabolic curve
exists.
\end{corol}


\begin{thebibliography}{99}
\bibliographystyle{plain}
\bibitem{ABT}M. Abate, F. Bracci and F. Tovena, \emph{Index theorems for holomorphic self maps}, to appear in  Ann. of Math.
%\bibitem{A}M.Abate, \emph{Diagonalization of nondiagonalizable discrete holomorphic dynamical systems},
% Amer.J. Math., 122, 2000,757-781.
\bibitem{AB}M. Abate, \emph{The residual index and the dynamics of holomorphic maps tangent to the identity},
Duke Math. J., (1) 107, 2001, 173-207.
\bibitem{BR}F. Bracci, \emph{The dynamics of holomorphic maps near curves of fixed points}, Ann. Sc. Norm. Sup. Pisa Cl. Sci. (5)
 vol. II (2003), 493-520.
\bibitem{BT}F. Bracci and F. Tovena, \emph{Residual indices of holomorphic maps relative to singular curves
of fixed points on surfaces}, Math. Z., 242, 2002, 481-490.
\bibitem{Bri}E.Brieskorn e H.Knorrer, \emph{Plane algebraic curves}, Birkhauser, 1986.
\bibitem{BRU}M. Brunella, \emph{Birational geometry of foliations}, First Latin American Congress
 of Math., IMPA, Rio de Janeiro, Brazil 2000.
\bibitem{CS}C. Camacho and P. Sad \emph{Invariant varieties through singularities of holomorphic vector fields},
Ann. of Math. (2) 115, 1982, 579-595.
\bibitem{CA}J. Cano, \emph{Construction of invariant curves for singular holomorfic vector fields}, Proc.
Amer. Math. Soc., 125(9) 1997,2649-2650.
%\bibitem{CG}l.Carleson e T.W.Gamelin, \emph{Complex Dynamics}, Springer, 1992.
%\bibitem{EC}J.\`Ecalle, \emph{Les fonctions r\'esurgentes, Tome III:l'equation du pont et la classification
%analytiques des objects locaux}, Publ. Math. Orsay, 85-5, 1985, Universit\'e de Paris-Sud, Orsay.
%\bibitem{GPV}G.Gentili, F.Podest\`a e E.Vesentini \emph{Lezioni di geometria differenziale}, Bollati
%Boringhieri, 1995.
%\bibitem{GH}P.A.Griffiths e J.Harris, \emph{Principles of algebraic geometry}, Wiley classics library,
% 16, 1981.
%\bibitem{GRI}P.A.Griffiths, \emph{Introduction to algebraic curves}, Traslations of mathematical monographs, 76, 1989.
\bibitem{GR}R. C. Gunning and H. Rossi, \emph{Analytic functions of several complex variables}, Prentice-Hall,Inc., 1965.
%\bibitem{HA}M.Hakim, \emph{Analytic transformations of $(\mathbb{C}^p,0)$ tangent to the identity},
% Duke. Math. J., 92,1998,403-428.
\bibitem{LAU}H. B. Laufer, \emph{Normal two dimensional singularities},  Ann.Math.Stud.71, 1971.
%\bibitem{L}D.Lehmann, \emph{Methods of differential geometry in analytic and algebraic geometry}, Note delle
%lezioni tenute presso l'Universit\`a di Roma ``Tor Vergata'' scritte da F.Bracci e G. Minervini, 2001,
%http://www.mat.uniroma2.it/~fbracci/web/downloads/lehmann.ps.
\bibitem{SW}T. Suwa, \emph{Indices  of holomorphic vector fields relative to invariant curves on
surfaces}, Proc. Amer. Math. Soc. 123, (1995), 2989-2997.
\bibitem{MM}J. F. Mattei and R. Moussu, \emph{Holonomie et int\'egrales premi\`eres}, Ann. Scient. de l'Ecole Norm. Sup\'erieure,
4 \`eme s\'erie-t.13, (1980), 469-523.
\bibitem{SB}M. Sebastiani, \emph{Sur l'existence de s\'eparatices locales des feuilletages des surfaces}, An.
 Acad. Bras.Ci., 69, 2, 1997, 159-162.
%\bibitem{SW}T.Suwa, \emph{Localization of characteristic classes and applications},
%Note delle lezioni tenute presso l'Universit\`a di Roma ``Tor Vergata'' scritte da F.Bracci, 2001,
% http://www.mat.uniroma2.it/~fbracci/web/downloads/suwa.ps.
%\bibitem{YANG}K.Yang,\emph{Complex algebraic geometry},Monographs and textbooks in pure and applied
% mathematics vol.149 , 1991.
\end{thebibliography}
\end{document}